\documentclass[12pt]{amsart}
\usepackage{amssymb}
\usepackage{amsfonts}
\usepackage{amsmath}
\usepackage{array}
\usepackage[ps,matrix]{xypic}

\usepackage{tableaux}

\def\QMR{{\it QMR}}
\def\MR{{\bf MR}}

\def\AKS{\operatorname{\mathcal H}}
\def\AKSz{\AKS_{n,r}(0)}
\def\AKSq{\AKS_{n,r}(q)}
\def\AKSzm{\AKS_{m,r}(0)}
\def\AKSdd{\AKS_{2,2}(0)}
\def\AKSqd{\AKS_{4,2}(0)}

\newtheorem{example}{Example}[section]

\newtheorem{theorem}[example]{Theorem}

\newtheorem{proposition}[example]{Proposition}

\newtheorem{lemma}[example]{Lemma}

\def\FQSym{{\bf FQSym}}

\def\NCSF{{\bf Sym}}
\def\QSym{{\it QSym}}

\def\<{\langle}
\def\>{\rangle}

\def\C{\operatorname{\mathbb C}}

\def\Z{\operatorname{\mathbb Z}}

\def\F{{\bf F}}

\def\G{{\bf G}}

\def\SG{{\mathfrak S}}

\def\Des{\operatorname{Des}}
\def\Des{C}

\title[Representation theory of the $0$-Ariki-Koike-Shoji algebras]
{Representation theory of the\\ $0$-Ariki-Koike-Shoji algebras}

\author[F. Hivert, J.-C.~Novelli, and J.-Y.~Thibon]%
{Florent Hivert, Jean-Christophe Novelli, and Jean-Yves Thibon}

\address[] {Institut Gaspard Monge, Universit\'e de Marne-la-Vall\'ee \\
5 Boulevard Descartes \\Champs-sur-Marne \\77454 Marne-la-Vall\'ee cedex 2 \\
FRANCE}
\email[Florent Hivert]{hivert@univ-mlv.fr}
\email[Jean-Christophe Novelli]{novelli@univ-mlv.fr}
\email[Jean-Yves Thibon]{jyt@univ-mlv.fr} 
\date{}

\begin{document}

\begin{abstract}
We investigate the representation theory of certain specializations of the
Ariki-Koike algebras, obtained by setting $q=0$ in a suitably normalized
version of Shoji's presentation.
We classify the simple and projective modules, and describe restrictions,
induction products, Cartan invariants and decomposition matrices.
This allows us to identify the Grothendieck rings of the towers of algebras in
terms of certain graded Hopf algebras known as the Mantaci-Reutenauer descent
algebras, and Poirier Quasi-symmetric functions.
\end{abstract}

\maketitle

\section{Introduction}

Given an \emph{inductive tower of algebras}, that is, a sequence of
algebras
\begin{equation}
A_0 \hookrightarrow A_1 \hookrightarrow \ldots \hookrightarrow A_n
\hookrightarrow \ldots,
\end{equation}
with embeddings $A_m\otimes A_n \hookrightarrow A_{m+n}$ satisfying an
appropriate associativity condition, one can introduce two \emph{Grothendieck
rings}
\begin{equation}
{\mathcal G}:=\bigoplus_{n\ge 0}G_0(A_n)\,,\quad
{\mathcal K}:=\bigoplus_{n\ge 0}K_0(A_n)\,,
\end{equation}
where $G_0(A)$ and $K_0(A)$ are the (complexified) Grothendieck groups of the
categories of finite-dimensional $A$-modules and projective $A$-modules
respectively, with multiplication of the classes of an $A_m$-module $M$ and an
$A_n$-module $N$ defined by
\begin{equation}
[M] \cdot [N] = [M\widehat{\otimes} N] =
[M\otimes N \uparrow_{A_m\otimes A_n}^{A_{m+n}}].
\end{equation}

On each of these Grothendieck rings, one can define a coproduct by means of
restriction of representations, turning these into mutually dual Hopf algebras.

The basic example of this situation is the character ring of symmetric groups
(over $\C$), due to Frobenius. Here the $A_n=\C\SG_n$ are semi-simple
algebras, so that
\begin{equation}
{G}_0(A_n) = {K}_0(A_n)= R(A_n),
\end{equation}
where $R(A)$ denotes the vector space spanned by isomorphism classes of
indecomposable modules which in this case are all simple and projective.
The irreducible representations $[\lambda]$ of $A_n$ are parametrized by
partitions $\lambda$ of $n$, and the Grothendieck ring is isomorphic to the
algebra $Sym$ of symmetric functions under the
correspondence
$[\lambda] \leftrightarrow s_\lambda$,
where $s_\lambda$ denotes the Schur function associated with $\lambda$.

Other known examples with towers of group algebras over the complex numbers,
$A_n=\C G_n$, include the cases of wreath products
$G_n = \Gamma\wr\SG_n$ (Specht), finite linear groups
$G_n = GL(n,\F_q)$ (Green), \emph{etc.}, all related to symmetric functions
(see~\cite{Mcd,Zel}).

Examples involving non-semisimple specializations of Hecke algebras have also
been worked out.
Finite Hecke algebras of type $A$ at roots of unity ($A_n=H_n(\zeta)$,
$\zeta^k=1$) yield quotients and subalgebras of $Sym$ \cite{LLT}
\begin{equation}
{\mathcal G} = Sym/(p_{km}=0), \quad
{\mathcal K} = \C\left[p_i\,|\,i\not\equiv0\pmod k\right],
\end{equation}
supporting level $1$ irreducible representations of the affine Lie algebra
$\widehat{sl}_k$, while Ariki-Koike algebras at roots of unity give rise to
higher level representations of the same Lie algebras \cite{Ari}.
The $0$-Hecke algebras $A_n=H_n(0)$ correspond to the pair Quasi-symmetric
functions/ Noncommutative symmetric functions, ${\mathcal G}=\QSym$,
${\mathcal K}=\NCSF$ \cite{NCSF4}.
Affine Hecke algebras at roots of unity lead to $U^+(\widehat{sl}_k)$ and
$U^+(\widehat{sl}_k)^*$ \cite{Ari}, and the cases of affine Hecke generic
algebras can be reduced to a subcategory admitting as Grothendieck rings
$U^+(\widehat{gl}_\infty)$ and $U^+(\widehat{gl}_\infty)^*$ \cite{Ari}.

A further interesting example is the tower of $0$-Hecke-Clifford algebras
\cite{Ols,BHT}, giving rise to the peak algebras \cite{NCSF2,Stem}.

Here, we shall show that appropriate versions at $q=0$ of the Ariki-Koike
algebras (presentation of Shoji \cite{Sho,Sasho}) admit as Grothendieck
rings two known combinatorial Hopf algebras, the Mantaci-Reutenauer descent
algebras (associated with the corresponding wreath products) \cite{MR},
and their duals, a generalization of quasi-symmetric functions, introduced
by Poirier in~\cite{Poi} and more recently considered in~\cite{NT-coul,BH}.

\bigskip
This article is structured as follows. We first define the
$0$-Ariki-Koike-Shoji algebras $\AKSz$, and introduce a special basis, well
suited for analyzing representations.
Next, we obtain the classification of simple $\AKSz$-modules, which turn out
to be all one-dimensional, and labelled by $r(r+1)^{n-1}$ combinatorial
objects called cyclotomic ribbons. We then describe induction products and
restrictions of these simple modules, which allows us to identify the first
Grothedieck ring $\mathcal G$ with a Hopf subalgebra of Poirier's
Quasi-symmetric functions, dual to the Mantaci-Reutenauer Hopf algebra.
This duality gives immediately the Grothendieck ring $\mathcal K$ associated
with projective modules. An alternative labelling of the indecomposable
projective modules leads then to a simple description of basic operations such
as induction products, restriction to $H_n(0)$, or induction from a
$H_n(0)$-projective module. Summarizing, we obtain an explicit description of
the Cartan-Brauer triangle, in particular of the Cartan invariants and of the
decomposition matrices.

In all the paper, we will make use of a set $C=\{1,\ldots,r\}$, called the
\emph{color set}.

\section{The $0$-Ariki-Koike-Shoji algebras}

In \cite{Sho}, Shoji obtained a new presentation of the Ariki-Koike algebras
defined in~\cite{AK}.
We shall first give a presentation very close to his, put $q=0$ in the
relations and then prove some simple results about another basis of the
resulting algebra. To get our presentation from the one of Shoji, one has to
replace $qa_i$ by $T_{i-1}$ and $q^2$ by $q$.

Let $u_1,\ldots,u_r$ be $r$ distinct complex numbers. We shall denote by
$P_k(X)$ the polynomial
\begin{equation}
\label{Lag}
P_k(X) := \prod_{1\leq l\leq r,\ l\not=k} \frac{X-u_l}{u_k-u_l}.
\end{equation}

Let $\AKSq$ be the associative algebra generated by the elements
$T_1, \ldots,T_{n-1}$ and $\xi_1,\ldots,\xi_n$ subject to the following
relations:

\begin{eqnarray}
(T_i-q)(T_i+1) = 0 \quad (1\leq i\leq n-1), \\
T_i T_{i+1} T_i = T_{i+1} T_i T_{i+1} \quad (1\leq i\leq n-2), \\
T_i T_j = T_j T_i \qquad (|i-j|\geq2), \\
{(\xi_j-u_1)\cdots (\xi_j-u_r)} =0 \qquad (1\leq j\leq n), \\
\xi_i \xi_j = \xi_j \xi_i \qquad (1\leq i,j \leq n), \\
T_i \xi_i = \xi_{i+1}T_i - (q-1)
\sum_{c_1<c_2}{(u_{c_2}-u_{c_1})P_{c_1}(\xi_i)P_{c_2}(\xi_{i+1})}
\quad (1\leq i\leq n-1), \\
T_i (\xi_{i+1}+\xi_i) = (\xi_{i+1}+\xi_i)T_i \quad (1\leq i\leq n-1)\\
T_i \xi_j = \xi_j T_i \qquad (j\not= i-1,i).
\end{eqnarray}

As noticed in~\cite{Sho}, from this presentation, it is obvious that a
generating set is given by the $\xi_1^{c_1}\cdots\xi_n^{c_n}\cdot T_\sigma$
with $\sigma\in\SG_n$ and $c_i$ such that $0\leq c_i\leq r-1$. Shoji proves
that this is indeed a basis of $\AKSq$.
Moreover, a simple adaptation of his proof enables us to conclude that this
property still holds a $q=0$.
This can also be directly proved thanks to the multiplication relations
between the Hecke generators and the Lagrange basis to be presented next.
This algebra $\AKSz$, which we call the $0$-Ariki-Koike-Shoji algebra, will be
our main concern in the sequel.

If ${\bf c}=(c_1,\ldots,c_n)$ is a word on $C$, we define
\begin{equation}
L_{\bf c} := P_{c_1}(\xi_1) \cdots P_{c_n}(\xi_n).
\end{equation}

Since the Lagrange polynomials (Equation~(\ref{Lag})) associated with $r$
distincts complex numbers are a basis of $\C_{r-1}[X]$ (polynomials of degree
at most $r-1$), the next proposition holds.

\begin{proposition}
The set 
\begin{equation}
\{ B_{{\bf c},\sigma} := L_{\bf c} T_\sigma\},
\end{equation}
where $\sigma\in\SG_n$ and ${\bf c}=(c_1,\ldots,c_n)$ is a color word, is a
basis of $\AKSz$.
\end{proposition}

Recall that a composition is any finite sequence of positive integers
$I=(i_1,\ldots,i_k)$. It can be pictured as a ribbon diagram, that is, a set
of rows composed of square cells of respective lengths $i_j$, the first
cell of each row being attached to the last cell of the previous one. $I$ is
called the \emph{shape} of his ribbon diagram.
Recall also that the descent composition $I=\Des(\sigma)$ of a permutation
$\sigma$ is the one whose diagram is obtained by writing the elements of
$\sigma$ one per cell so that the rows are weakly increasing and the columns
are strictly decreasing (French notation).

We shall represent the basis element $B_{{\bf c},\sigma}$ by a filling of a
composition diagram as follows: the composition is $\Des(\sigma)$ and its
$i$-th cell is filled with $c_i$ and $\sigma_i$.

Let us now describe the product by a generator on the left of a basis element:
on $B_{{\bf c},\sigma}$, the generator $\xi_i$ acts diagonally by
multiplication by $c_i$, so that it only remains to explicit the product of
$T_i$ by $L_{\bf c}$. One finds

\begin{equation}
T_i L_{\bf c} = L_{{\bf c}\sigma_i} T_i +
\left\{
\begin{array}{lcl}
- L_{{\bf c}} & \text{if} &c_i<c_{i+1}, \\
0 & \text{if} &c_i=c_{i+1}, \\
L_{{\bf c}\sigma_i} & \text{if} &c_i>c_{i+1}. \\
\end{array}
\right.
\label{tl-lt}
\end{equation}
where $\sigma_i$ acts on the right of $\bf c$ by exchanging $c_i$ and
$c_{i+1}$.

In fact the previous expression is the specialization $q=0$ of the following
apparently unnoticed relation in $\AKSq$:
\begin{equation}
T_i L_{\bf c} = L_{{\bf c}\sigma_i} T_i - (q-1)\left\{
\begin{array}{lcl}
-L_{{\bf c}} & \text{if} &c_i<c_{i+1}, \\
0 & \text{if} &c_i=c_{i+1}, \\
L_{{\bf c}\sigma_i} & \text{if} &c_i>c_{i+1}. \\
\end{array}
\right.
\end{equation}

This description enables us to analyze the left regular representation of
$\AKSz$ in terms of our basis elements. As a first application, we shall
obtain a classification of the simple $\AKSz$-modules.

\section{Simple modules}

\subsection{Definition}

Let $I$ be a composition of $n$ and ${\bf c}\in C^n$ be a color word of length
$n$. The pair $[I,{\bf c}]$ will be called a \emph{colored ribbon}, and
depicted as the filling of $I$ whose natural row reading is $\bf c$.
We say that this filling is a \emph{cyclotomic ribbon} (\emph{cycloribbon}
for short) if it is weakly \emph{increasing} in rows and weakly
\emph{decreasing} in columns.
Notice that there are $r(r+1)^{n-1}$ cycloribbons since when building the
ribbon cell by cell, one has $r$ possibilities for its first cell and then
$r+1$ possibilities for the next ones: $1$ possibility for the $r-1$ different
choices than the previous one and $2$ possibilities for the same choice
(either go right or go down). Here are the five cycloribbons of shape $(2,1)$
with two colors.

\begin{equation}
\label{5cyclo}
\smalltableau{1&1\\\ &1}
\qquad
\smalltableau{1&2\\\ &1}
\qquad
\smalltableau{1&2\\\ &2}
\qquad
\smalltableau{2&2\\\ &1}
\qquad
\smalltableau{2&2\\\ &2}
\end{equation}

\subsection{Eigenvalues}

Let $I$ be a composition of $n$ and ${\bf c}=(c_1,\ldots,c_n)$ be a color
word.
We say that their associated ribbon is a \emph{anticyclotomic ribbon}
(\emph{anticycloribbon} for short) if it is weakly \emph{decreasing} in rows
and weakly \emph{increasing} in columns. There are as many anticycloribbons
as cycloribbons. The relevant bijection $\phi$ from one set to the other is
the restriction of an involution on all colored ribbons: read a ribbon $R$ by
rows from top to bottom (French notation), and build the corresponding ribbon
cell by cell in the following way:
\begin{itemize}
\item if the $(i+1)$-th cell has the same content as the $i$-th cell, glue it
in the same position as in $R$ (right or down),
\item if the $(i+1)$-th cell does not have the same content as the $i$-th
cell, glue it in the other position (right or down).
\end{itemize}

For example, the colored ribbons below are exchanged by application of $\phi$:

\begin{figure}[ht]
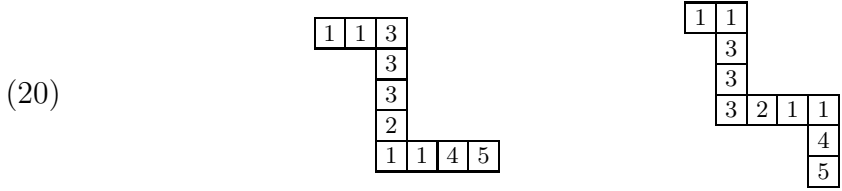

\begin{equation}
\smalltableau{1&1&3\\\ &\ & 3\\\ &\ &3\\\ &\ &2\\\ &\ &1&1&4&5}
\qquad \qquad \qquad
\smalltableau{1&1\\\ &3\\\ &3\\\ &3&2&1&1\\\ &\ &\ &\ &4\\\ &\ &\ &\ &5}
\end{equation}
\caption{A cycloribbon $R$ and its corresponding anticycloribbon
$\phi(R)$.}
\end{figure}

Let $[I,{\bf c}]$ be a cycloribbon.
Let $\eta_{[I,{\bf c}]}$ be the element of $\AKSz$ defined as
\begin{equation}
\eta_{[I,{\bf c}]} := L_{\bf c}\ \eta_I,
\end{equation}
where $\eta_I$ is the generator of the simple $H_n(0)$-module associated with
$I$ in the notation of \cite{NCSF4} Proposition~5.3 (see also~\cite{Nor}).
We will show later that this element generates a simple $\AKSz$-module.

\begin{proposition}
Let $[I,{\bf c}]$ be a cycloribbon.
Then
\begin{equation}
\eta_{[I,{\bf c}]} = \eta_{I'} L_{\bf c'},
\end{equation}
where $[I',{\bf c'}]$ are the shape and the mirror image of the reading of the
anticycloribbon $\phi([I,{\bf c}])$.
\end{proposition}

\begin{theorem}
\label{simples}
Let $[I,{\bf c}]$ be a cycloribbon. Then
\begin{equation}
S_{[I,{\bf c}]} := \AKSz \eta_{[I,{\bf c}]}
\end{equation}
is a simple module of $\AKSz$ realized as a minimal left ideal in its left
regular representation. The eigenvalue of $\xi_i$ is $c_i$, and that of $T_i$
is $-1$ or $0$ according to whether $i$ is a descent of the shape of
$\phi(R')$ or not.
All these simple modules are pairwise non isomorphic and of dimension $1$.
Moreover, all simple $\AKSz$-modules are one-dimensional and isomorphic to
some $S_{[I,{\bf c}]}$.
\end{theorem}

\begin{proof}
We only give a sketch of the proof.

The first part of the theorem follows directly from the previous propositions.
Next, we prove that any simple module of dimension $1$ is isomorphic
to some $S_{[I,{\bf c}]}$. This comes easily from Equation~(\ref{tl-lt}).

It remains to show that all simple modules are of dimension $1$. This is
done by the same argument as in~\cite{BHT} in our context: compute the
composition factors of the modules induced from simple modules of the
$0$-Hecke algebra. This is described in the next subparagraph.
\end{proof}

\subsection{Induction of the simple $0$-Hecke modules}

To describe the induction process, we will need a partial order on the
fillings of ribbons. Let $I$ be a composition and ${\bf c}=(c_1,\ldots,c_n)$.
The covering relation of the order $\leq_I$ corresponds to sort in increasing
order any two adjacent elements in the rows of $I$ or to sort in decreasing
order any two adjacent elements in the columns of $I$.
For example, the elements smaller than
\begin{equation}
T := \smalltableau{2&1\\\ &1&3&3\\\ &\ &\ &4\\\ &\ &\ &3}
\end{equation}
are
\begin{equation}
\smalltableau{1&2\\\ &1&3&3\\\ &\ &\ &4\\\ &\ &\ &3}
\qquad \qquad
\smalltableau{2&1\\\ &1&3&4\\\ &\ &\ &3\\\ &\ &\ &3}
\qquad \qquad 
\smalltableau{1&2\\\ &1&3&4\\\ &\ &\ &3\\\ &\ &\ &3}
\end{equation}

If $I$ is a composition of $n$, let $S_I:=H_n(0)\eta_I$ be the corresponding
simple module of $H_n(0)$ (notation as in~\cite{NCSF6}) and
\begin{equation}
M_I := S_I \uparrow_{H_n(0)}^{\AKSz}.
\end{equation}

Clearly, $M_I$ has dimension $r^n$ and admits $L_{\bf c} \eta_I$ as linear
basis, when $\bf c$ runs over color words.
For ${\bf c}\in C^n$, let $M_{{\bf c},I}$ be the $\AKSz$-submodule of $M_I$
generated by $L_{\bf c}\eta_I$.

\begin{lemma}
\begin{equation}
M_{{\bf c},I} \subseteq M_{{\bf c'},I} \quad\Longleftrightarrow\quad
{\bf c}\leq_I {\bf c'}.
\end{equation}
\end{lemma}

\begin{proof}
Let $i\in\{1,\ldots,n-1\}$.

\begin{itemize}
\item If $c_i<c_{i+1}$, and $T_i$ acts by $0$ on $\eta_I$, we get
$(1+T_i)L_{\bf c}\eta_I = L_{{\bf c}\sigma_i}T_i\eta_I=0$.
\item If $c_i<c_{i+1}$, and $T_i$ acts by $-1$ on $\eta_I$, we get
$(1+T_i)L_{\bf c}\eta_I = L_{{\bf c}\sigma_i}T_i\eta_I
=- L_{{\bf c}\sigma_i}\eta_I$, and so $-(1+T_i)$ sorts in decreasing order in
columns $c_i$ and $c_{i+1}$.
\item If $c_i=c_{i+1}$, then $T_iL_{\bf c}\eta_I = L_{\bf c}T_i\eta_I$ so
the result is either $0$ or $-1$ times $L_{\bf c}\eta_I$.
\item If $c_i>c_{i+1}$, and $T_i$ acts by $-1$ on $\eta_I$, we get
$T_iL_{\bf c}\eta_I = L_{{\bf c}\sigma_i}(1+T_i)\eta_I = 0$.
\item If $c_i>c_{i+1}$, and $T_i$ acts by $0$ on $\eta_I$, we get
$T_iL_{\bf c}\eta_I = L_{{\bf c}\sigma_i}(1+T_i)\eta_I
= L_{{\bf c}\sigma_i}\eta_I$, and so $T_i$ sorts in increasing order in rows
$c_i$ and $c_{i+1}$.
\end{itemize}
\end{proof}

The $M_{{\bf c},I}$ such that the pair ${\bf c},I$ is a cycloribbon are
simple.
It follows easily from the previous lemma that these form the socle of
$M_I$, and that $M_I$ admits a composition series involving only
one-dimensional modules. Since any simple $\AKSz$-module must appear as a
composition factor of some $M_I$, this proves Theorem~\ref{simples}.

\subsection{First Grothendieck ring}

Let us first recall that the composition factors of the induction product of
two simple $0$-Hecke modules labelled by compositions $I$ and $J$ are easily
described by means of permutations: let $\sigma$ (resp. $\tau$) be the
permutation with the maximum inversion number in the descent class $I$
(resp. $J$). Then the $0$-Hecke graph associated with both simple modules is
the same as the graph of the weak order restricted to the shifted
shuffle of the inverses of $\sigma$ and $\tau$. The composition
factors are then the descent compositions of the elements of this shifted
shuffle product.

In other words, the composition factors of the induction product of
two simple $0$-Hecke modules can be obtained by computing a product in
$\FQSym$ (compute a shifted shuffle of permutations) and taking the image of
each element by the morphism that sends a permutation to its descent
composition (same morphism as taking the commutative image of the usual
realization as noncommutative polynomials).
Since this is an epimorphism onto $\QSym$, this proves that the direct sum
of all Grothendieck groups associated with the $0$-Hecke algebra, is
isomorphic to $\QSym$.

\medskip
The case at hand can be worked out in a similar way. Instead of permutations,
one has to consider colored permutations as in~\cite{NT-coul}. Given a
cycloribbon $[I,{\bf c}]$, one associates with it the unique colored
permutation $(\sigma,u)$ where $\sigma$ is the permutation with the maximum
inversion number of the descent class $I$ and $u={\bf c}$. Recall that colored
permutations index the Hopf algebra $\FQSym^{(r)}$ \cite{NT-coul} and that
the product of the $\F$ basis (dual to the $\G$ basis, considering the colors
as a group) is given by the shifted shuffle of colored permutations.

Let $[I,{\bf c}]$ and $[J,{\bf c'}]$ be two cycloribbons.
Then the graph associated with the induction product of the corresponding
simple modules is the same as the graph of the shifted shuffle of the inverses
(considering the colors as a cyclic group) of the colored permutations
associated with both cycloribbons. Both graphs have the same structure, the
edges being given by $T_i$ or $1+T_i$ depending on whether $c_i\geq c_{i+1}$
or not (see Equation~(\ref{tl-lt})).

The simple module corresponding to a given colored permutation $(\sigma,u)$ is
the cycloribbon associated with the colored descent composition of
$(\sigma,u)$. Recall that the colored descent composition of a colored
permutation is the pair $(I,v)$ where $I$ is the unique composition which
descents are either a descent of $\sigma$ or a change of color in $u$, and
$v$ is the color of all elements of the corresponding block of the $k$-th row
of~$I$.

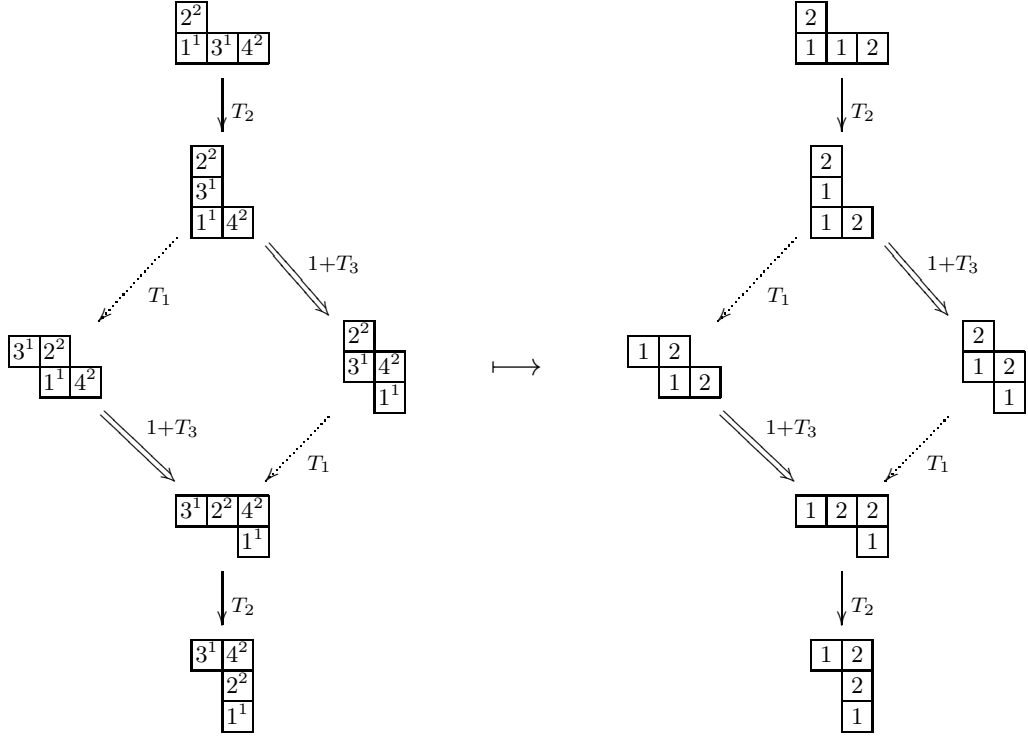
\begin{figure}[ht]
\newdimen\vcadre\vcadre=0.2cm 
\newdimen\hcadre\hcadre=0.2cm 
\def\GrTeXBox#1{\vbox{\vskip\vcadre\hbox{\hskip\hcadre%
      $\smalltableau{#1}$%
   \hskip\hcadre}\vskip\vcadre}}

$\vcenter{\xymatrix@R=0.7cm@C=0.6cm{
   & *{\GrTeXBox{2^2\\ 1^1 & 3^1 & 4^2}} \ar @{->}[d]^{T_2} & \\
   & *{\GrTeXBox{2^2\\ 3^1\\ 1^1 & 4^2}}
   \ar @{..>}[dl]^{T_1} \ar @{=>}[dr]^{1+T_3}& \\
 *{\GrTeXBox{3^1& 2^2\\\ &1^1 & 4^2}}\ar @{=>}[dr]^{1+T_3}
 & & *{\GrTeXBox{2^2\\ 3^1 & 4^2\\\ & 1^1}} \ar @{..>}[dl]^{T_1}\\
   & *{\GrTeXBox{3^1& 2^2& 4^2\\\ &\ &1^1}} \ar @{->}[d]^{T_2}& \\
   & *{\GrTeXBox{3^1& 4^2\\\ & 2^2\\\ &1^1}} & \\
}}
\qquad\longmapsto\qquad
\vcenter{\xymatrix@R=0.7cm@C=0.6cm{
   & *{\GrTeXBox{2\\ 1 & 1 & 2}} \ar @{->}[d]^{T_2} & \\
   & *{\GrTeXBox{2\\ 1\\ 1 & 2}}
   \ar @{..>}[dl]^{T_1} \ar @{=>}[dr]^{1+T_3}& \\
 *{\GrTeXBox{1& 2\\\ &1 & 2}}\ar @{=>}[dr]^{1+T_3}
 & & *{\GrTeXBox{2\\ 1 & 2\\\ & 1}} \ar @{..>}[dl]^{T_1}\\
   & *{\GrTeXBox{1& 2& 2\\\ &\ &1}} \ar @{->}[d]^{T_2}& \\
   & *{\GrTeXBox{1& 2\\\ & 2\\\ &1}} & \\
}}
$
\caption{\label{ex-induc}Induction product of two simple modules of $\AKSdd$
to $\AKSqd$}
\end{figure}

For example, Figure~\ref{ex-induc} presents the induction product of the two
simple modules $[(11),21]$ and $[(2),12]$ of $\AKSdd$ to $\AKSqd$.
In particular, we get the following composition factors written as
cycloribbons: $[(1,3),2112]$, $[(1,1,2),2112]$, $[(2,2),1212]$,
$[(1,2,1),2121]$, $[(3,1),1221]$, $[(2,1,1),1221]$. Notice that this
construction gives an effective algorithm to compute the composition factors
of the induction product of two simple modules.

\subsection{The Quasi-symmetric Mantaci-Reutenauer algebra}

Let $\QSym^{(r)}$ be the algebra of level $r$ quasi-symmetric functions,
defined in~\cite{Poi} (see also~\cite{NT-coul}).
This algebra is indexed by \emph{signed compositions}, (or \emph{vector
compositions} in~\cite{NT-coul}) which we will not
define. Poirier remarks (Lemma~11, p. 325) that the functions indexed by
\emph{descent signed compositions} (a subset of the previous set) can be
computed through the operation "take the commutative image" of a
noncommutative series, very similar to the definition of
$\FQSym$~\cite{NCSF6}.
It happens that these functions (indexed by descent signed compositions) form
a subalgebra $\QMR^{(r)}$ of $\QSym^{(r)}$ (see~\cite{BH}). We obtained
this algebra \cite{NT-coul} as a quotient of $\QSym^{(r)}$ by killing some
monomials.
We will prove that the Grothendieck ring of the tower of $\AKSz$ algebras is
isomorphic to $\QMR^{(r)}$.

There is a trivial bijection between descent signed compositions and
cycloribbons, so that we can speak of $F_{[I,{\bf c}]}$ without ambiguity.

\begin{lemma}[\cite{BH,NT-coul}]
The $F_{[I,{\bf c}]}$, where $[I,{\bf c}]$ runs over the cycloribbons,
span a Hopf subalgebra $\QMR^{(r)}$ of $\QSym^{(r)}$, isomorphic to the dual
of the Mantaci-Reutenauer algebra $\MR^{(r)}$ defined in~\cite{MR}.
\end{lemma}

Let 
\begin{equation}
{\mathcal G}^{(r)} := \bigoplus_{n\geq 0} G_0(\AKSz)
\end{equation}
be the Grothendieck ring of the tower of $\AKSz$ algebras. Define a
characteristic map
\begin{equation}
\begin{array}{rcl}
ch: & {\mathcal G}^{(r)} &\to \QMR^{(r)}\\
& [S_{[I,{\bf c}]}] &\mapsto F_{[I,{\bf c}]}.
\end{array}
\end{equation}

\begin{theorem}
$ch$ is an isomorphism of Hopf algebras.
\end{theorem}

\begin{proof}
Both ${\mathcal G}^{(r)}$ and $\QMR^{(r)}$ are the image of $\FQSym^{(r)}$
over the morphism that sends any permutation to its colored descent
composition: for ${\mathcal G}^{(r)}$, this follows from the previous
paragraph, whereas for $\FQSym^{(r)}$, it is proved in~\cite{NT-coul}.
\end{proof}

\section{Projective modules}

The previous results immediately imply, by duality, a description of the
Grothendieck ring of the category of projective $\AKSz$-modules.
As a Hopf algebra,
\begin{equation}
{\mathcal K} = \bigoplus_{n\geq0} K_0(\AKSz) \simeq \MR^{(r)}
\end{equation}
is isomorphic to the Mantaci-Reutenauer algebra, and under this isomorphism,
the class of indecomposable projective modules are mapped to the subfamily of
the dual basis $F_{[I,{\bf c}]}$ of Poirier quasi-symmetric functions
labelled by colored descent sets, or cycloribbons.

This labelling is however not always convenient and it is useful to introduce
another one, by \emph{colored compositions} $(I,u)$, that is, pairs formed by
a composition $I=(i_1,\ldots,i_p)$, and a color word $u=(u_1,\ldots,u_p)$ of
the same length.

The bijection between colored compositions and anticycloribbons is easy to
describe: starting with a colored composition, one rebuilds an
anticycloribbon by separating adjacent blocks of the colored composition with
different colors and gluing them back in the only possible way according to
the criterion of being an anticycloribbon. Conversely, one separates two
adjacent blocks of different colors and glue them one below the other.

Recall that the $\MR^{(r)}$ Hopf algebra can be defined (see~\cite{NT-coul})
as the free associative algebra over the symbols
$(S_{j}^{(i)})_{j\geq1 ; 1\leq i\leq r}$, graded by $\deg S_j^{(i)}=j$, and
with coproduct
\begin{equation}
\Delta S_n^{(k)} = \sum_{i=0}^n S_i^{(k)} \otimes S_j^{(k)}.
\end{equation}
For a colored composition $(I,u)$ as above, one defines
\begin{equation}
S^{(I,u)} := S_{i_1}^{(u_1)} \cdots S_{i_p}^{(u_p)}.
\end{equation}
Clearly, the $S^{(I,u)}$ form a linear basis of $\MR^{(r)}$.
There is a natural order on colored compositions, which generalizes the
anti-refinement order on ordinary compositions: one says that
$(J,v)\leq (I,u)$ if $(J,v)$ can be obtained from $(I,u)$ by adding up groups
of consecutive parts of the same color.

For example, with two colors, the anti-refinements of
$(2\overline{1}\overline{2}13)$ are

\begin{equation}
(2\overline{1}\overline{2}13), 
(2\overline{3}13), 
(2\overline{1}\overline{2}4), 
(2\overline{3}4).
\end{equation}

The \emph{colored ribbon basis} $R_{(I,u)}$ of $\MR^{(r)}$ can now be defined
by the condition

\begin{equation}
S^{(I,u)} =: \sum_{(J,v)\leq (I,u)} R_{(J,v)}.
\end{equation}

Let $(I,u) = (i_1,\ldots, i_p ; u_1,\ldots, u_p)$ and
$(J,v)= (j_1,\ldots, j_q;v_1,\ldots, v_q)$ be two colored compositions.
We set
\begin{equation}
(I,u) \cdot (J,v) := (I\cdot J,u\cdot v),
\end{equation}
where $a\cdot b$ denotes the concatenation.

Moreover, if $u_p=v_1$, we set
\begin{equation}
(I,u)\triangleright (J,v) :=
(i_1,\ldots, i_{p-1},(i_{p}+j_1),j_2,\ldots,j_q ;
 u_1,\ldots, u_p,v_2,\ldots, v_q).
\end{equation}

The colored ribbons satisfy the very simple multiplication rule:
\begin{equation}
R_{(I,u)} R_{(J,v)} = R_{(I,u)\cdot (J,v)} +
\left\{
\begin{array}{ll}
R_{(I,v)\triangleright(J,v)} & \text{if}\ u_p=v_1,\\
0 & \text{if}\ u_p\not=v_1.
\end{array}
\right.
\end{equation}

Let $[K,{\bf c}]$ be an anticycloribbon. Let $P_{[K,{\bf c}]}$ be the
indecomposable projective module whose unique simple quotient is the simple
module labelled by $\phi([K,{\bf c}])$ and let $(I,u)$ be the corresponding
colored composition.

\begin{theorem}
\label{K2MR}
The map
\begin{equation}
\begin{array}{rcl}
{\bf ch}: & {\mathcal K}     &\to \MR^{(r)} \\
          &[P_{[K,{\bf c}]}] &\mapsto R_{(I,u)}
\end{array}
\end{equation}
is an isomorphism of Hopf algebras.
\end{theorem}

The main interest of the labelling by colored compositions is that it allows
immediate reading of some important information. For example, it follows from
Theorem~\ref{K2MR} that the products of complete functions $S^{(I,u)}$ are the
characteristics of the projective $\AKSz$-modules obtained as induction
products in which each factor $S_m^{(i)}$ is the characteristic of the
one-dimensional projective $\AKSzm$-module on which all the $T_i$ act by $0$
and all the $\xi_j$ by the same eigenvalue $u_i$.

We see that, as in the case of $H_n(0)$, each indecomposable projective
$\AKSz$-module occurs as a direct summand of such an induced module, and that
the direct sum decomposition is given by the anti-refinement order.
For example, the identity
\begin{equation}
S^{(2\overline{1}\overline{2}13)} =
R_{ (2\overline{1}\overline{2}13)} +
R_{(2\overline{3}13)} +
R_{(2\overline{1}\overline{2}4)} +
R_{(2\overline{3}4)}
\end{equation}
indicates which indecomposable projective direct summands compose the
projective $\AKS_{9,2}(0)$-module defined as the outer tensor product
\begin{equation}
S_2 \widehat\otimes S_{\overline{1}} \widehat\otimes S_{\overline{2}}
\widehat\otimes S_1 \widehat\otimes S_3.
\end{equation}

Their restrictions to $H_9(0)$ (and hence, their dimensions) can be computed
by means of the following result:

\begin{theorem}
The homomorphism of Hopf algebras
\begin{equation}
\begin{array}{rcl}
\pi: &\MR^{(r)}    & \to \NCSF \\
     & S_{j}^{(i)} & \mapsto S_j
\end{array}
\end{equation}
maps the class of a projective $\AKSz$-module to the class of its restriction
to $H_n(0)$.
\end{theorem}

Continuing the previous example, we see that the restriction of
$P_{2\overline{1}\overline{2}13}$ to $H_9(0)$ is given by
\begin{equation}
\pi(R_{2\overline{1}\overline{2}13}) =
\pi(R_2R_{\overline{12}}R_{13}) =
R_2R_{12} R_{13} =
R_{21213} + R_{2133} + R_{3213} + R_{333}.
\end{equation}

Dually, one can describe the induction of projective $H_n(0)$-modules
to $\AKSz$.

\begin{theorem}
Let $I$ be a composition of $n$ and let $N_I$ be the $\AKSz$-module induced by
the indecomposable projective $H_n(0)$-module $P_I$.
Then
\begin{equation}
N_I \simeq \bigoplus P_{[I,{\bf c}]},
\end{equation}
where the sum runs over all the anticycloribbons of shape $I$.
\end{theorem}

For example, let us complete the case $I=(2,1)$ with two colors. The
following five anticycloribbons appear in the induction of $P_I$, with
respective dimensions $3$, $6$, $3$, $2$, and $2$
\begin{equation}
\smalltableau{2&1\\\ &1}
\qquad
\smalltableau{2&1\\\ &2}
\qquad
\smalltableau{1&1\\\ &2}
\qquad
\smalltableau{2&2\\\ &2}
\qquad
\smalltableau{1&1\\\ &1}
\end{equation}

\medskip
Finally, we can describe the Cartan invariants and the decomposition matrices
by means of the maps

\begin{equation}
\begin{array}{rcl}
e: &\MR^{(r)} &\to Sym^{(r)} = (Sym)^{\otimes r} \simeq Sym(X_1,\ldots,X_r),\\
   & S_j^{(i)} &\mapsto h_j(X_i)
\end{array}
\end{equation}
and
\begin{equation}
\begin{array}{rcl}
d:& Sym^{(r)} &\hookrightarrow \QMR^{(r)}\\
& h_j(X_i) &\mapsto F_{[(j),i^j]}.
\end{array}
\end{equation}

Then, the Cartan map is $c=d\circ e$, and the entry
$C_{[I,{\bf c}], [J,{\bf d}]}$ of the Cartan matrix (giving the multiplicity
of the simple module $S_{[J,{\bf d}]}$ as a composition factor of the
indecomposable projective module $P_{[I,{\bf c}]}$) is equal to the
coefficient of $F_{[J,{\bf d}]}$ in $c(R_{[I,{\bf c}]})$.
The decomposition map is given by $d$ and the decomposition matrix expresses
the tensor product of Schur functions
$S_\lambda = s_{\lambda^{(1)}}(X_1) \cdots s_{\lambda^{(r)}}(X_r)$
on the Poirier basis $F_{[I,{\bf c}]}$.


%
%
%
%
%
%
%
%
%
%

\end{document}